\documentclass[11pt]{amsart}
\usepackage{amssymb}
\usepackage{graphicx} 
\usepackage{enumerate}
\usepackage{mathtools}
\usepackage{color,soul}
\usepackage{cite}
\usepackage{tikz}
\usetikzlibrary{matrix}
\usepackage{caption}
\usepackage{amsmath, amsthm}
\usepackage[colorlinks=true, allcolors=blue]{hyperref}
\usepackage[noabbrev]{cleveref}
 \usepackage{relsize}
 \usetikzlibrary{cd}
 \usetikzlibrary{decorations.pathreplacing}
 \usepackage{tikz,calc}
\usepackage{color}
\usepackage[margin=1.2in]{geometry}

\usepackage{multicol,todonotes}

\newcommand{\theoremname}{Theorem:}

\title{Danceability, Directed by Braid Index}
\author{Sol Addison} 
\author[]{Nancy Scherich}
\author[]{Lila Snodgrass}
\date{}					

\begin{document}

\begin{abstract}

Schaffer introduced the concept of \emph{danceability} of a knot diagram. In this paper, we expand upon Schaffer's ideas to create a danceability knot invariant and show that this invariant is bounded above by the braid index. 

\end{abstract}

\maketitle

\thispagestyle{empty}

\section*{Knots and Danceability} 

\begin{figure}[htp]
\centering
\begin{picture}(220,70)
\put(-100,0){\includegraphics[scale=.3]{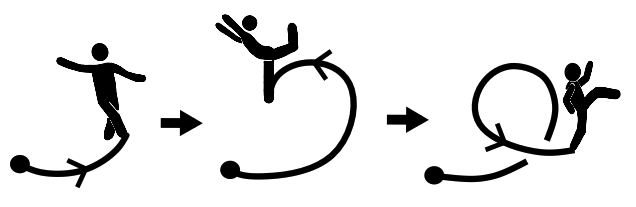}\hspace{.6cm}}
\put(140,0){\includegraphics[scale=.3]{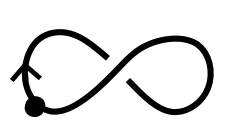}\hspace{.7cm}}
\put(260,0){\includegraphics[scale=.3]{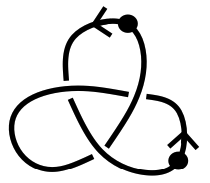}\hspace{.7cm}}
\put(-120,5){(a)}
\put(110,5){(b)}
\put(230,5){(c)}
\end{picture}
\caption{(a) Depiction of a dancer tracing a knot following the under-first rule. (b) A 1-danceable
diagram. (c) The Trefoil knot diagram is 2-danceable.}  \label{fig:1}
\end{figure}
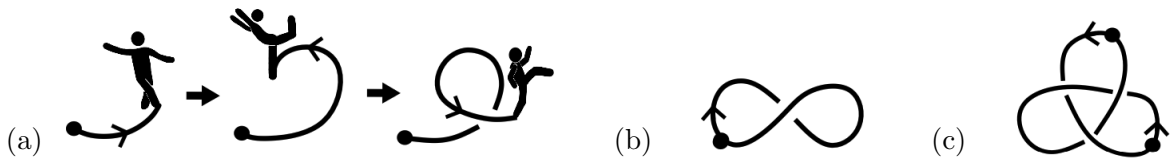

In Mathematics, a \emph{knot} is a closed curve in 3 dimensional space with no self intersections. We can think of knots as being made of extremely thin and flexible material. If you bend, stretch, or rearrange a knot without cutting it, we would say  it is the same knot, it just looks different. 
A \emph{knot diagram} is a planar projection of the knot with no triple points and where the over and under strands of a crossing are denoted in the diagram, as shown in Figure \ref{fig:1} (c).

In his 2021 Bridges article \cite{KS}, Karl Schaffer introduced the concept of \emph{danceability}, asking which knot diagrams can be danced. That is, a knot is \emph{danceable}  if you can choose an orientation, or directional flow, of the knot and a point on the knot so that  a dancer can start at the point, dance in the directional flow of the knot, and traverse the entire knot with the restriction that the dancer must pass through every crossing as the under-strand first. 
We will refer to this restriction as the ``under-first" rule.
Visually, you can imagine the dancer's path drawing the knot over time, and later times must always travel on top of earlier times, as shown in Figure \ref{fig:1} (a). 
If a point can be found on the diagram satisfying these conditions, the diagram is said to be \emph{solo-danceable}, or \emph{1-danceable}.

Schaffer showed that any knot diagram that is 1-danceable is \emph{unknotted}, that is the knot can be untied, without cutting it, to have no crossings. 
He gave a  very intuitive proof using the visual of a ribbon and a maypole. More technically, he described a Morse lifting of the knot diagram to a 3 dimensional closed curve which has two critical values, and thus must be unknotted \cite[Theorem 10]{Schar}. 
All nontrivial knot diagrams are not 1-danceable and must, therefore, be at least duet- or 2-danceable. 
Schaffer explained that for a knot diagram to be \emph{2-danceable}, both dancers must start independently at different points, move along the path in the same direction, and end at the other dancer’s initial point, while still following the under-first rule. Naturally, one might ask what happens when a knot diagram cannot be traversed by just two dancers? 
To pick up where Schaffer left off, we formalize Schaffer’s definition of $n$-danceability in such a way that we can define, and study, a danceability knot invariant.

\vspace{.25cm}

\noindent \textbf{Definition 1.}\emph{
For $n\geq 2$, an oriented knot diagram is \emph{n-danceable} if there exists $n$ initial points on the knot diagram so that the following conditions are true.}

\begin{enumerate}
    \item \emph{One dancer starts at each of the $n$ initial points, for a total of $n$ dancers.}

    \item  \emph{Each dancer travels in the direction of the pre-chosen orientation of the diagram, and stops dancing when they reach the next initial point. This ensures that the total path traveled by all dancers traverses the entire knot diagram with no overlap.}


    \item \emph{The speed that each dancer travels can vary and can be chosen so that the simultaneous tracing of the diagram by all $n$ dancers follows the under-first rule at every crossing.}
\end{enumerate}

Since two different diagrams for the same knot can look wildly different, as explained by Reidemeister's Theorem \cite{R}, we want to step away from diagram-dependency. 
Thus, we define the \emph{danceability} of a knot $K$, denoted $da(K)$, to be the minimal number of dancers, $n$, for which there is a diagram of $K$ that is $n$-danceable. It is clear from the minimality in the definition that danceability is a knot invariant. That is, for any two knots $K_1$ and $K_2$, if $da(K_1)\neq da(K_2)$ then $K_1$ and $K_2$ must be distinct knots.
\vspace{.25cm}

\noindent \textbf{Theorem 1.}\emph{ Danceability is a knot invariant.}

\section*{Danceability is Bounded Above by Braid Index}

Given a knot, the danceability is the minimum number of dancers required to dance the knot. There are other knot invariants that also use minimality, such as the crossing number. The \emph{crossing number} of a knot is the  minimal number of crossings in any knot diagram for the knot.
As Schaffer pointed out, any knot diagram with $n$-crossings is $n$-danceable by simply assigning a dancer to the underpass of each crossing. This shows that danceability is bounded above by the crossing number. However, danceability is not equal to crossing number as the Trefoil has crossing number 3, but is 2-danceable, as shown in Figure \ref{fig:1} (c).
\vspace{.25cm}

\noindent \textbf{Theorem 2.} [Schaffer \cite{KS}] \emph{ Danceability is bounded above by the crossing number.}
\vspace{.25cm}

\begin{figure}[htp]
\centering
\begin{picture}(220,70)
\put(-60,0){\includegraphics[scale=.3]{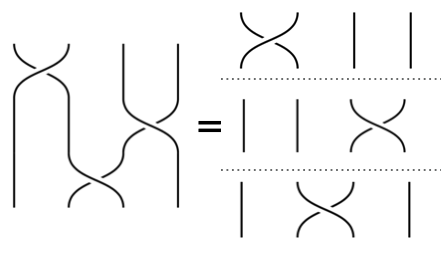}\hspace{.6cm}}
\put(130,0){\includegraphics[scale=.4]{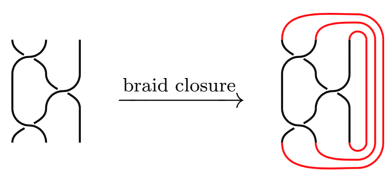}\hspace{.7cm}}
\put(-80,5){(a)}
\put(110,5){(b)}
\end{picture}
\caption{(a) A braid shown as a stack of 1-crossing braids. (b) The closure of a braid.}  \label{fig:braids}
\end{figure}

 To further investigate the minimality of danceability, we can explore braids. 
 A \emph{braid on $n$ strands} is a collection of $n$ arcs in $\mathbb{R}^3$ that connect $n$ top points with height 1 to $n$ bottom points with height 0 in such a way that all arcs monotonically flow up. See Figure \ref{fig:braids} for an example. 
 The monotonicity requirement implies that each crossing in the braid can be arranged to occur at a different height. Moreover,  every braid can be seen as a vertical stack of smaller braids, each with one crossing.
 We can take the \emph{closure} of a braid to get a knot, or a link, by identifying the top and bottom end points with new arcs, see Figure \ref{fig:braids} (b). 
 A crucial theorem in braid theory is Alexander's Theorem which states that every knot can be represented as the closure of a braid \cite{A}.

 Furthermore, for any knot $K$, the \emph{braid index}, denoted $b(K)$, is the fewest number of strands needed to express $K$ as a closed braid \cite{ARTIN}. Similar to the crossing number, we find that the braid index of a knot serves as an upper bound for the knot’s danceability. 


\vspace{.25cm}

\noindent \textbf{Theorem 3.}
    \emph{Danceability is bounded above by braid index. That is, for any knot $K$, $da(K)\leq b(K)$.}

\vspace{.25cm}

\begin{proof}

 Suppose $K$ is a knot diagram that is the closure of a braid, $\beta$, with $n\geq 2$ strands. Suppose also that $b(K)=n$ and $\beta$ realizes the minimal braid index of $K$. Choose the orientation of $K$ induced by the upwards orientation of the braid.
 Place $n$ dancers, one on each strand at the bottom of the braid. All dancers travel “up” the braid.
After the dancers traverse the braid, they travel the arcs in the closure, which connect back down to the initial starting points. This shows conditions (1) and (2) in the definition of $n$-danceable are satisfied. See Figure \ref{fig:braid_Examples} (a).

\begin{figure}[htp]
\centering
\begin{picture}(220,120)
\put(130,30){\includegraphics[scale=.19]{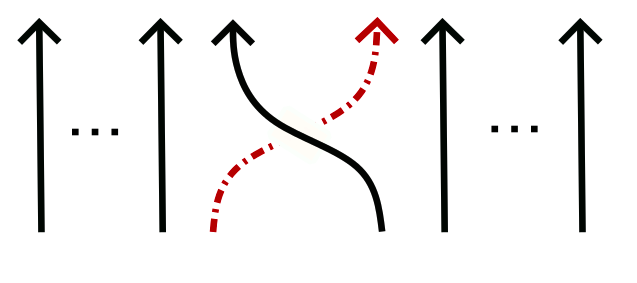}\hspace{.6cm}}
\put(-40,0){\includegraphics[scale=.18]{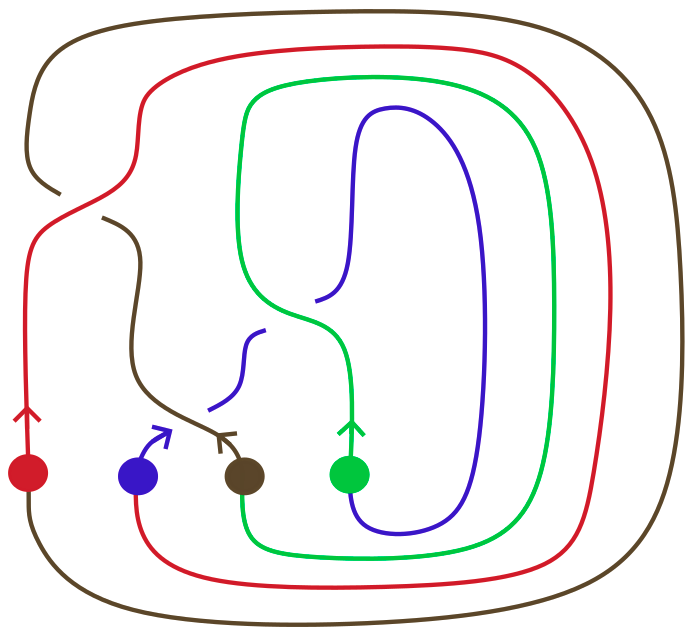}\hspace{.6cm}}
\put(-80,5){(a)}
\put(110,5){(b)}
\end{picture}
\caption{(a) Starting positions and orientation for a braid closure.  (b) The dancer on the under-strand, shown in dashed red, can travel faster than the dancer on the over-strand, satisfying under-first. }  \label{fig:braid_Examples}
\end{figure}

An $n$-stranded braid with a single crossing can be danced with $n$ dancers by having all dancers traverse upwards at the same speed, except for the dancer on the under-strand of the crossing. 
That dancer speeds up to pass under the crossing, before the dancer on the over-strand arrives, and then slows down again to end the braid at the same time as the remaining dancers. This shows that a single crossing braid can be danced to satisfy the under-first rule. See Figure \ref{fig:braid_Examples} (b).

We can arrange the crossings of the $\beta$ to occur one crossing at a time in the height direction, and view $\beta$ as a stack of $n$-stranded braids with a single crossing. Since the dancers can traverse one crossing at a time and satisfy the under-first rule, we have shown condition (3) in the definition of $n$-danceable is satisfied.


 Thus, $K$’s danceability (minimum number of dancers required) is bounded above by $b(k)$. 
\end{proof}

 Is the inequality in Theorem 3 strict; is $da(K)<b(K)$ for all $K$? No! For example, the Trefoil knot $T$, as shown in Figure \ref{fig:1}, is 2-danceable and is not the unknot, so $da(T)=2$.  $T$ is known to have braid index 2, which shows that $da(T)=b(T)$. Since all nontrivial knots must have danceability at least 2, we get the following corollary.


\vspace{.25cm}

\noindent \textbf{Corollary 1.}
    \emph{Any knot $K$ with $b(K)=2$ must also have $da(K)=2$.}

\begin{proof}
If $K$ has braid index 2, then $K$ is not the unknot, so $2\leq da(K)$. By Theorem 3, we see that 
$2\leq da(K)\leq b(K)=2.$ 
\end{proof}

\begin{figure}[h]
    \centering
    \includegraphics[scale=.3]{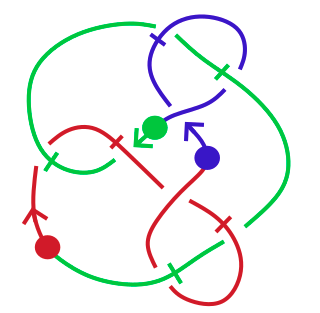}
    \caption{The $8_{15}$ knot has $b(8_{15})=4$, $br(8_{15})=3$, and $da(8_{15})=3$. The three dancer's starting points are shown as solid dots. The paths are marked with a small perpendicular line segment at points on the path where a dancer may need to wait for another dancer to pass under before continuing. }
    \label{fig:8_15}
\end{figure}

For example, all $T(2,q)$ torus knots have braid index 2 and therefore must have danceability  2. Is danceability always equal to the braid index? No! For example, the  $8_{15}$ knot shown in Figure \ref{fig:8_15} is 3-danceable, so $da(8_{15})\leq 3$. However, $8_{15}$ is known to have braid index 4, so $da(8_{15})<b(8_{15})$.


\section*{Future Directions}

 A knot diagram is \emph{alternating} if the crossings alternate under, over, under, over, as one travels along the knot. Since a dancer must satisfy the under-first rule, when traveling along an alternating knot, the dancer would have to wait at every over-crossing for another dancer to dance under the crossing first. We believe this means a sufficiently complicated alternating knot will require a large number of dancers.  

\vspace{.25cm}

\noindent \textbf{Conjecture 1.}
   \emph{If $K$ is an alternating knot with $3 \leq b(K)$, then $3\leq da(K)\leq b(K)$.}

\vspace{.25cm}

Since danceability is defined as a minimum over all possible diagrams of a knot, it is a difficult invariant to compute. It would be very useful to find a lower bound for the danceability. Through computational trial and error, we suspect that the bridge index could be a lower bound.  The \emph{bridge index} of a knot is the minimum over all knot diagrams of the number of local maxima of the knot, where the knot is viewed as a smooth closed curve in 3 dimensional space. For example, the $8_{15}$ knot has bridge index 3 and braid index 4. We showed in Figure \ref{fig:8_15} that $8_{15}$ is 3-danceable, and we suspect that it cannot be 2-danceable since the bridge index is 3. A referee for this article pointed out a potential argument that the bridge index could be an upper bound for the danceability of a knot. Let the bridge index of a knot $K$ be denoted by $br(K)$. It is well known that  $br(K)\leq b(K)$, and yet from our current understanding, it is unclear where danceability fits into this inequality.











\bigskip

\section*{Acknowledgements} This project is part of an undergraduate student research experience at Elon University lead by the second author.
We would like to thank our anonymous referees for their helpful comments.


\bibliography{Danceability}
\bibliographystyle{plain}

\end{document}